\documentclass[twoside,12pt,reqno]{amsart}
\usepackage{amsmath,amsthm, amsfonts,amssymb}
\usepackage{braket}
\usepackage{color,soul,mathtools}
\usepackage[dvipsnames]{xcolor}
\usepackage[all]{xy}
\usepackage{graphicx}
\usepackage{indentfirst}
\usepackage{bm}
\usepackage{mathrsfs}
\usepackage{latexsym}
\usepackage{hyperref}

\newtheorem{Thm}{Theorem}[section]

\newtheorem{corollary}[Thm]{\bf Corollary}

\newtheorem{theorem}[Thm]{\bf Theorem}

\title[Positivity under Semidefinite linear Constraints]{Certain Linear Combinations of Exponential Functions are Positive under Semidefinite linear constraints}
\author{
	Robert~Lin
}
\address{	Department of Physics \\
	Harvard University \\
	Cambridge, MA 02138}
\date{\today}

\begin{document}

\maketitle
\begin{abstract}
	In this article, I introduce a group-theoretical method to prove positivity of certain linear combinations (with coefficients generally lying in $\mathbb{C}$) of exponential functions under a set of semidefinite linear constraints. The basic group-theoretic fact we rely on is the positivity of the fusion coefficients for multiplication of group characters.
\end{abstract}

\section{Introduction}
In prior work joint with Jonathan Boretsky\cite{LB}, it was shown that the positivity of certain linear combinations (with coefficients generally lying in $\mathbb{C}$) of exponential functions (parameterized by a single variable $t$) implies a set of semi-definite linear constraints. In this work, I will show that the converse is also true, and that in fact, fundamentally, the reason for the positivity of these combinations is due to \textit{group theory}.

In addition to demonstrating an unexpected connection between semidefinite optimization and group theory, the method I use to prove positivity may be of general interest as well. Possible applications include the study of inequalities arising in the study of semigroups, such as hypercontractive bounds.

\section{The Main Theorem (Case where $G=\mathbb{Z}_N$)}
Suppose $P_t$ is a semigroup that acts on the group algebra $\mathcal{L}G$ by $P_t \ket{g} = e^{-t |g|} \ket{g}$, where $g \in G$ and $|\cdot|:G\rightarrow \mathbb{R}_{\geq 0}$. Assume that $|e|=0$ and $|g|=|g^{-1}|$ for all $g \in G$. In \cite{LB}, the sum representation  $P_t = \sum_{i=1}^N p_i(t) \sigma_i$ was introduced, where $\sigma_i$'s are multipliers in the group element basis induced by the corresponding irreducible characters $\chi_i$. 

In this section, for simplicity, I will present my main theorem in the case $G=\mathbb{Z}_N$. Thus, $\sigma_i$ is defined by $\sigma_i \ket{j} = \chi_{ij} \ket{j}$, where $\chi_{ij} = q^{ij}$ and $q= \exp(2\pi i /N)$. The $\sigma_i$ form a set of mutually orthogonal vectors in the Hilbert space $B(\mathcal{L}G)$ defined by $\langle a, b \rangle = \text{Tr}(a^* b)$. 

The choice of a length function $|\cdot|$ on the group defines a set of functions $p_i(t)$. We assume that such a decomposition exists for all $t\geq 0$, which is a constraint on the semigroup $P_t$. For $\mathbb{Z}_N$, the existence of a decomposition is guaranteed \cite{LB} since $|g|=|g^{-1}|$ by assumption. 

We can now state the theorem:
\begin{theorem}
	\label{maintheorem}
	If there exists $\epsilon >0$ such that for all $i \in \mathbb{Z}_N$, $p_i(t)\geq 0$ for all $t\leq \epsilon$, then for any $j \in \mathbb{Z}_N$, $p_j(t) \geq 0$ for all $t\geq 0$.
\end{theorem}
\begin{proof}
	For any $t>0$, take $n$ large such that $t/n< \epsilon$. Then, $P_t = P_{n(\frac{t}{n})} = \left(P_{\frac{t}{n}}\right)^n = \left( \sum_{i=1}^{N} p_i(t/n) \sigma_i\right)^n$ since $P_t$ is a semigroup. Since $\sigma_{a_1} \sigma_{a_2} \cdots \sigma_{a_n} = \sigma_{a_1+ a_2 + \cdots + a_n}$, we can replace the product by 
	\begin{equation}
		\sum_{a_1,\cdots, a_n=0}^{N-1} p_{a_1}(t/n) p_{a_2}(t/n) \cdots p_{a_n}(t/n) \sigma_{a_1+\cdots + a_n}.
	\end{equation}
	The coefficients in the product are all \textbf{nonnegative} by hypothesis. Since summing over products of nonnegative numbers always yields nonnegative numbers, we get by summing over all indices $a_j$ for fixed $a_1+\cdots + a_n=i$ that $p_i(t) \geq 0$ for all $t\geq 0$. 
\end{proof}
\begin{corollary}
	It suffices to check for all $i$ the values of $p_i(t)$ and $p_i'(t)$ at $t=0$ in order to check whether $p_j(t)\geq 0$ for all $t\geq 0$. Since $p_j(t)$ is a linear combination of exponentials of the form $e^{-|g|t}$, it follows that we only need to check a set of semidefinite constraints to show that $p_j(t)$ are all positive.
\end{corollary}

\textbf{Example: }For $\mathbb{Z}_6$, $p_{i\neq 0}=0$ and $p_{0}(0)=0$. The inequalities forced by $p'_{i\neq 0} (t=0) \geq 0$ are clearly necessary. The main theorem implies that these are sufficient. To be explicit, we represent the length function on $\mathbb{Z}_6$ by setting $|0|=0$, $|1|=|5|=a$, $|2|=|4|=b$, and $|3|=c$. Since the length function completely specifies the semigroup $P_t$, which has the sum representation $P_t = \sum_{i=0}^{5} p_i(t) \chi_i$, we can solve for the $p_i$'s, yielding
\begin{align}
	p_1(t) &= \frac{1}{6} \left(1+e^{-at} - e^{-bt} - e^{-ct} \right) \\
	p_2(t) &= \frac{1}{6} \left(1 - e^{-at} - e^{-bt} + e^{-ct} \right) \\
	p_3(t) &= \frac{1}{6} \left(1- 2 e^{-at} + 2 e^{-bt} - e^{-ct} \right).
\end{align}
We omit $p_0(t)$ which is a sum of positive quantities \cite{LB} (since $\chi_0$ corresponds to the trivial representation), hence always positive. 

Application of Theorem \ref{maintheorem} tells us that if for all $j \in \mathbb{Z}_N$, $p_j'(0)\geq 0$, then for all $i \in \mathbb{Z}_N$, $p_i(t) \geq 0$ for all $t\geq0$.
For $a,b,c\geq 0$,  there are only three inequalities resulting from $p_{i\neq 0}'(0)\geq 0$ (since $p_i = p_{6-i}$):
\begin{align}
	a- b- c &\leq 0 \\
	-a - b + c &\leq 0\\
	-2a + 2b - c &\leq 0.
\end{align}
Thus, $p_i(t)$ are all nonnegative for $a,b,c$ satisfying the above inequalities.

This set of inequalities for $p_i(t)$ is quite nontrivial, since generally speaking, one cannot deduce $p_i(t) \geq 0$ simply from knowledge that $p_i'(0) \geq 0$. 

\section{Generalization to Arbitrary Finite Groups}
The above result can be generalized to any finite group. Let $\chi_r$ be the character for an irrep labeled by $r$, and for each irrep $r$, define $\sigma_r \ket{g} = \chi_r(g) \ket{g}$. Further assume that $|\cdot|$ is a class function. Then $P_t$ admits a decomposition as $P_t = \sum_r p_r(t) \sigma_r$.

\begin{theorem}
	\label{maintheorem2}
	If there exists $\epsilon >0$ such that for all irreps $r$, $p_r(t)\geq 0$ for all $t\leq \epsilon$, then for any irrep $s$, $p_s(t) \geq 0$ for all $t\geq 0$.
\end{theorem}
\begin{proof}
		For any $t>0$, take $n$ large such that $t/n< \epsilon$. Then, $P_t = P_{n(\frac{t}{n})} = \left(P_{\frac{t}{n}}\right)^n = \left( \sum_{i \text{ irrep}} p_i(t/n) \sigma_i\right)^n$ since $P_t$ is a semigroup. Positivity follows now since the tensor product of irreducible representations of a finite group can be completely reduced, and the multiplicity of the irreducible representation $\rho_c$ in $\rho_a \otimes \rho_b$ is precisely $n_{ab}^c$, which is always nonnegative. By extension, we can write $\sigma_{a_1} \sigma_{a_2} \cdots \sigma_{a_n} = n_{a_1 a_2 \ldots a_n}^{b} \sigma_b$, where $n_{a_1 a_2 \ldots a_n}^{b}$ is the multiplicity of the irrep $\rho_b$ in $\rho_{a_1} \otimes \rho_{a_2} \otimes \cdots \rho_{a_n}$.  Thus, the product can be rewritten as 
\begin{equation}
	\sum_{a_1,\cdots, a_n} \sum_{b} p_{a_1}(t/n) p_{a_2}(t/n) \cdots p_{a_n}(t/n) n_{a_1 a_2 \ldots a_n}^{b} \sigma_{b}.
\end{equation}
	The coefficients in the product are all \textbf{nonnegative} by hypothesis. Since summing over products of nonnegative numbers always yields nonnegative numbers, we get by summing over all indices $a_j$ for fixed $b$ that for any irrep s, $p_s(t) \geq 0$ for all $t\geq 0$. 
\end{proof}

\section{Acknowledgments}

I wish to thank Jonathan Boretsky for helpful review of my write-up of my result.

I acknowledge the research support of ARO grant W911NF-20-1-0082 through the MURI project ``Toward Mathematical Intelligence and Certifiable Automated Reasoning: From Theoretical Foundations to Experimental Realization."

\end{document}